\documentclass[10 pt, conference]{article}

\usepackage{graphicx}      
\usepackage[round]{natbib}

\usepackage{diagbox}

\usepackage{url}

\usepackage[htt]{hyphenat}
\usepackage{amssymb}
\usepackage{bm}
\usepackage{siunitx}
\usepackage[nosumlimits]{amsmath}
\usepackage{color}

\begin{document}
\title{System Identification for Hybrid Systems using Neural Networks}

\author{Mattias F{\"a}lt, Pontus Giselsson\thanks{The authors are with the Department of Automatic Control, Lund University, Sweden. Both authors are financially supported by Swedish Research Council. E-mail: \texttt{\{mattiasf,pontusg\}@control.lth.se}.
This work has been submitted to IFAC for possible publication.}
}

\maketitle
\begin{abstract}
With new advances in machine learning and in particular powerful learning libraries, we illustrate some of the new possibilities they enable in terms of nonlinear system identification. For a large class of hybrid systems, we explain how these tools allow for identification of complex dynamics using neural networks. We illustrate the method by examining the performance on a quad-rotor example.
\end{abstract}


\section{Introduction}
In this paper we study the possibility of identifying dynamics of a hybrid system using a neural network.
We present how modern tools can be used to achieve this and
and illustrate the performance on a quad-rotor example.

Identification of continuous-time system dynamics has several advantages to discretized alternatives.
The dynamics of physical systems are often inherently driven by continuous dynamics,
making modeling and identification natural in continuous time.
Not only does this mean that parameters in a continuous model closely relate to real physical properties,
but also that the inherent structure, such as sparsity, can be captured in a continuous time model,
which could otherwise be lost in a discretized version.

However when parts of a model is unknown,
a parameterization of the uncertainty
is needed for identification.
Neural networks are well known to have the universal function approximator property, \cite{cybenko-1989, hornik-1991}
which have made them ubiquitous.

Identifying dynamics as part of an \emph{ordinary differential equation} (ODE)
using neural using neural networks was done already in the early 1990s, e.g. \cite{martinez-1992}.
The approach was known as Runge-Kutta neural networks, where a neural network was cloned and connected in a way where
the output of the connection corresponded to the solution to the ODE.
This form of the connection enables explicit formulas for the back-propagation of gradients,
allowing the neural network to be trained.

However, ODEs are not always enough to accurately model dynamics.
Differential algebraic equations are needed to correctly model many systems, for example in robotics.
And in many cases, discrete events are affecting the dynamics, requiring hybrid models.
For example in many real settings, identification in open loop is not possible,
and a discrete time controller has to be considered in the model.

These problems would be hard to handle with the Runge-Kutta neural network approach.
Even if some hybrid dynamics could be captured with a well chosen step-size,
where the timings of discrete events are known a priori,
deriving the back-propagation rules would be tedious and problem dependent.

Recent interest in using ODEs in neural networks, see \cite{chen-2018}, has generated several tools for automatic differentiation of differential equation solvers.
These tools it is possible to simulate, and back-propagate gradients through a variety of specialized
solvers, including methods with adaptive step-sizes and for stiff problems, DEAs, DDEs and many more. 
In this paper we show how these tools can be used to identify dynamics with neural networks even for complicated models.
To illustrate this we study an example of a quad-rotor in closed loop, with a discrete time controller, where we identify a drag model.

\section{Method}
Let the form of a general hybrid system be given by the dynamics
\begin{align}
    \dot{\bm{x}}(t) &= f(\bm{x}(t),t,\bm{\theta}) &  \textrm{if } (\bm{x}(t), t) \not\in D \nonumber\\
    \bm{x}(t^+) &= h(\bm{x}(t^-),t^-,\bm{\theta}) &  \textrm{if } (\bm{x}(t), t) \in D \label{eq:hybrid-system}\\
    \bm{x}(0) &= \bm{x_0}(\bm{\theta})&\nonumber
\end{align}
where $f$ represents the continuous dynamics, $h$ the discrete dynamics triggered by the guard $D$, $\bm{x}(t^-)$ and $\bm{x}(t^+)$ are the states just before and after a discrete update $h$, and $\bm{\theta}$ are a set of parameters to be estimated.
In this work we consider the guard to be of the form $D=\cup D_i$, where each $D_i$ is of one of the forms
\begin{align*}
   D_i&=\{(\bm{x},t) \mid g_i(\bm{x},t)=0\}           && \textrm{where }g_i\textrm{ continuous}\\
   D_i&=\{(\bm{x},t) \mid g_i(\bm{x},t)=0, t\in T_i\} && \textrm{where }T_i\textrm{ finite}.
\end{align*}

This form of hybrid systems is captures a wide range of hybrid dynamics. The first form of guard captures conditions triggered by the evolution of the state. The bouncing ball is a classical example of this, which could be captured by the guard function $g_i(h)=h$ where $h$ is the height of the ball w.r.t ground.
The second case captures any type of condition where the timings of events are known.
This can for example capture the dynamics of a discrete-time control system as in the example in Section \ref{sec:example}.

The parameter $\bm{\theta}$ is not restricted to represent quantities in the model.
It could for example be used to define the layers in a neural network,
and therefore serve for general function approximation.
It is be possible to model the whole $f$ function as a neural network and identify the dynamics this way.
However, when parts of the dynamics is known, it is naturally better to include this knowledge into the model, and identify only the unknown parts.
We also allow the initial state to depend on $\bm{\theta}$, this way the framework captures both the case where the initial state is known and when it is unknown.

We assume a setting where some loss function $L$ is given so that the best parameters $\bm{\theta^*}$, in some sense, are given by
\begin{align*}
\bm{\theta}^* = \underset{\bm{\theta}}{\textrm{argmin}}\, L(\bm{\theta}),
\end{align*}
where $L$ is a function of the states $\bm{x}$,
that depend on $\bm{\theta}$ via the dynamics \eqref{eq:hybrid-system}.
$L$ depends also usually on several other variables, such as initial values, measurements, and reference signals.
We do not explicitly denote these dependencies.
Solving this problem analytically is not possible in general. We will rely on some type of gradient descent method. For this we need the derivatives $\partial L/\partial \bm{\theta}$.
Since L depends on $\bm{\theta}$ via the states $\bm{x}$, we need to be able to differentiate through the hybrid dynamics \eqref{eq:hybrid-system} to find these derivatives.

For continuous-time systems, implementations for automatic differentiation of combinations of neural networks and ODE solvers now exist in several languages.
A common approach is to solve the adjoint equations to calculate the sensitivity $\partial \bm{x}(t)/\partial \bm{\theta}$.
Although there has been research on sensitivity analysis using adjoint methods for hybrid systems, e.g. \cite{Zhang-2017},
we are not aware of any convenient software that implements this approach.

The Julia language see \cite{Julia-2017}, has a wide range of feature rich and efficient differential equation solvers that natively allow for automatic differentiation.
The \texttt{DifferentialEquations.jl} package by \cite{Rackauckas-2017} contains over 140 native ODE solvers, and special solvers for Split ODE, Second order ODE, SDE, DDE, DAE, BVP, Jump Problems,
most or all of which allow for automatic differentiation using either forward- or reverse-mode automatic differentiation.
We use this package, together with the machine learning library \textrm{Flux.jl} to train our model, see \cite{innes-2018}.
We use the package \textrm{DiffEqFlux.jl} by \cite{Rackauckas-2019} to combine these tools in a convenient manner.

The callback event handling interface in \texttt{Differential-Equations.jl}
allows for handling of the discrete dynamics triggered by the guards $D_i$. The \texttt{ContinuousCallback} handles cases where $g_i$ are continuous functions,
where the crossings of $0$ are found by a root-finding method.
The \texttt{DiscreteCallback} can be used to handle events triggered at specific times.

\section{Example}\label{sec:example}
To illustrate the capabilities of this framework, we study a quad-rotor example.
The goal is to identify a set of parameters as well as a drag model,
based on data generated in closed-loop with a discrete time zero-order hold controller.

Identification in closed loop is inherently problematic due to lack of excitation, see \cite{ljung-1999}.
To alleviate this,
we excite the system with an extra excitation signal on the inputs.
We assume that we are only able to measure the positions of the quad-rotor at specific time points,
and that these measurements are noisy.
To simplify the example,
we use a state feedback controller and assume that the internal controller has access to the true state.
%
\subsection{System}
We use the simple quad-rotor model given in \cite{sabatino-2015}, with parameters chosen based on the very the small Crazyflie quad-rotor in \cite{Greiff-2017}.
The dynamics are
\begin{align*}
\dot{\eta_1}= & \xi_1 + \eta_3\cos(\eta_1)\tan(\eta_2)+\xi_2\sin(\eta_1)\tan(\eta_2)\\
\dot{\eta_2}= & \xi_2\cos(\eta_1)-\xi_3\sin(\eta_1)\\
\dot{\eta_3}= & \xi_3\cos(\eta_1)/\cos(\eta_2)+\xi_2\sin(\eta_1)/\cos(\eta_2)\\
\dot{\xi_1}= & \xi_3\xi_2(I_{y}-I_{z})/I_{x}+\tau_1/I_{x}\\
\dot{\xi_2}= & \xi_1\xi_3(I_{z}-I_{x})/I_{y}+\tau_2/I_{y}\\
\dot{\xi_3}= & \xi_2\xi_1(I_{x}-I_{y})/I_{z}+\tau_3/I_{z}\\
\dot{v_1}= & \xi_3v_2-\xi_2v_3-g\sin(\eta_2)+f_{w_1}/m\\
\dot{v_2}= & \xi_1v_3-\xi_3v_1+g\sin(\eta_1)\cos(\eta_2)+f_{w_2}/m\\
\dot{v_3}= & \xi_2v_1-\xi_1v_2+g\cos(\eta_2)\cos(\eta_1)-f_{t}/m+f_{w_3}/m\\
\dot{p_1}= & R_{1}(\bm{\eta})\bm{v}\\
\dot{p_2}= & R_{2}(\bm{\eta})\bm{v}\\
\dot{p_3}= & R_{3}(\bm{\eta})\bm{v}
\end{align*}
where $R$ is the rotation matrix from the coordinate system defined by the arms of the quad-rotor (local frame) to that of the room (global frame):
\begin{align*}
  R(\bm{\eta})=\begin{bmatrix}R_{1}(\bm{\eta})\\
    R_{2}(\bm{\eta})\\
    R_{3}(\bm{\eta})
  \end{bmatrix}=
  \begin{bmatrix}
  c_{2}c_{3} &
  c_{3}s_{1}s_{2}-c_{1}s_{3} &
  s_{1}s_{3}+c_{1}c_{3}s_{2}\\
  c_{2}s_{3} &
  c_{1}c_{3}+s_{1}s_{2}s_{3} &
  c_{1}s_{3}s_{2}-c_{3}s_{1}\\
  -s_{2}
  & c_{2}s_{1}
  & c_{1}c_{2}
  \end{bmatrix},
\end{align*}
where $s_i := \sin(\eta_i)$ and $c_i := \cos(\eta_i)$.
An overview of the variables and parameters is presented in Tables \ref{tab:vars} and \ref{tab:params}.
The drag is modeled as
\begin{align*}
f_{w_1} & = d_1(v_1-(R(\bm{\eta})^T\bm{w}))_1)\\
f_{w_2} & = d_2(v_2-(R(\bm{\eta})^T\bm{w}))_2)\\
f_{w_3} & = d_3(v_3-(R(\bm{\eta})^T\bm{w}))_3)
\end{align*}
where $\bm{w}$ is wind in the global reference frame, $d_i(x)=-k_i x|x|$ where $k_i$ are \emph{drag parameters}.
The thrust and torques are given by the motor dynamics
\begin{align*}
  f_t & = b(\Omega_1^2 + \Omega_2^2 + \Omega_3^2 + \Omega_4^2)\\
  \tau_1 & = bl(\Omega_3^2 - \Omega_1^2)\\
  \tau_2 & = bl(\Omega_4^2 - \Omega_2^2)\\
  \tau_3 & = d(\Omega_2^2 + \Omega_4^2 - \Omega_1^2 - \Omega_3^2)\\
\end{align*}
where $\bm{\Omega}=(\Omega_1,\Omega_2,\Omega_3,\Omega_4)$ are the control signals.
To simplify notation somewhat, we denote the full state vector $\bm{x}:=(\bm{\eta}, \bm{\xi}, \bm{v}, \bm{p})$.
\begin{table}
   \caption{\label{tab:vars}Variables in the quad-rotor model.}
   \centering
   \begin{tabular}{| c c c c |}
      \hline
      Description & Variable & Frame & Unit \\
      \hline\hline
      State & $\bm{x}=(\bm{\eta}, \bm{\xi}, \bm{v}, \bm{p})$ & - & - \\
      \hline
      Euler angles & $\bm{\eta}=(\eta_1,\eta_2,\eta_3)$ & Local & \si{\radian} \\
      \hline
      Ang. velocity & $\bm{\xi}=(\xi_1,\xi_2,\xi_3)$ & Local & \si{\radian\per\second} \\
      \hline
      Velocity & $\bm{v}=(v_1,v_2,v_3)$ & Local & \si{\metre\per\second} \\
      \hline
      Position & $\bm{p}=(p_1,p_2,p_3)$ & Global & \si{\metre} \\
      \hline
      Torque & $\bm{\tau}=(\tau_1,\tau_2,\tau_3)$ & Local & \si{\newton\metre} \\
      \hline
      Drag & $\bm{f_w}=(f_{w_1},f_{w_2},f_{w_3})$ & Local & \si{\newton} \\
      \hline
      Rotor & $\bm{\Omega}=(\Omega_1,\Omega_2,\Omega_3,\Omega_4)$ & - & \si{\per\minute} \\
      \hline
      Wind & $\bm{w}=(w_1,w_2,w_3)$ & Global & \si{\metre\per\second} \\
      \hline
      Inertia & $\bm{I}=(I_1,I_2,I_3)$ & Local & \si{\kilogram\per\square\metre}$\cdot10^{-6}$ \\
      \hline
      Drag param. & $\bm{k}=(k_1,k_2,k_3)$  & - & - \\
      \hline
      Thrust & $f_t$ & Local & \si{\newton} \\
      \hline
   \end{tabular}
   \vspace{0.5cm}
   \caption{\label{tab:params}Parameters in the quad-rotor model.}
   \begin{tabular}{| c | c | c |}
      \hline
      Parameter & Value & Unit \\
      \hline\hline
      $\bm{w}$ & $(0.223, 0.354, -0.154)$ & \si{\metre\per\second} \\
      \hline
      $\bm{I}$ & $(6.48,6.48,9.98)$  & \si{\kilogram\per\square\metre}$\cdot10^{-6}$ \\
      \hline
      $\bm{k}$ & $(0.03618,0.03618,0.1809)$ & - \\
      \hline
      $b,d,l$ & $(2.2\cdot10^{-8},10^{-9},0.046)$ & (-\si{},-\si{},\si{\metre}) \\
      \hline
      $m,g$ & $(0.027,9.81)$ &  (\si{\kilogram},\si{\metre\per\square\second}) \\
      \hline
   \end{tabular}
\end{table}

\subsubsection{Controller}
A simple affine, discrete-time state-feedback controller is used to stabilize the process.
Let $T = \{t_k\}_{k=1}^N$ be the set of sample times.
For the reminder of the paper, we will use indexing with brackets to represent times in a discrete-time signal,
and parentheses for continuous-time variables.
The control signal is then defined as
\begin{align*}
\bm{\Omega}(t) ¤= \bm{\Omega}[k] \quad \forall t \in [t_k, t_{k+1})
\end{align*}
with
\begin{align}\label{eq:update-control}
   \bm{\Omega}[k] = \bm{\Omega}_0 + K(\bm{x}(t_k) - \bm{r}[k]) + \bm{r_n}[k]
\end{align}
where $\bm{\Omega}_0$ is an offset, $K\in\mathbb{R}^{4\times 12}$,
$\bm{r}[k]$ is the reference signal, and $\bm{r_n}[k]$ is a known excitation signal of known
deterministic noise.

To cast this into the hybrid system form \eqref{eq:hybrid-system} we add $\bm{\Omega}$ to the state vector $\bm{x}$ with $\dot{\bm{\Omega}}=\bm{0}$, define the guard
\begin{equation*}
D = \{(x,t) \mid t\in T\}
\end{equation*}
and set the update equation $h$ according to \eqref{eq:update-control},
which can be implemented with the \texttt{DiscreteCallback} functionality.

\subsection{Problem Formulation}
We assume a setting where a number of trajectories of the closed-loop system is collected,
and the goal is to identify the dynamics.
We let $T=\{0,0.05, 0.1,\dots, 10\}$ and generate $M=15$ trajectories for training,
with the process and controller given in the previous section.
We now assume that we don't know the inertia parameters $\bm{I}$ or
the drag forces $\bm{f_w}$, and set the goal to identify these.

The data used for training are noisy measurement of the positions
\begin{align*}
   \bm{y}[k] = \bm{p}(t_k) + \bm{n}[k]
\end{align*}
where $\bm{n}[k]\in\mathbb{R}^3$ is zero-mean Gaussian noise,
as well as the reference signals $\bm{r}[k]$, the excitation signal $\bm{r_n}[k]$ and the initial positions $\bm{x}(0)$.
For full details of the generation of these signals, see Appendix \ref{app:data}.

Let $\hat{\bm{y}}_i[k]$ be the estimated outputs of trajectory $i$, given some guess $\bm{\theta}$ of the inertia and drag model.
We set the loss function $L$ to be the least squares cost over all trajectories
\begin{align}\label{eq:cost-prediction}
L(\bm{\theta}) = \frac{1}{MN}\sum_{i=1}^M\sum_{k=1}^N \|\hat{\bm{y}}_i[k]-\bm{y}_i[k]\|_2^2.
\end{align}
so that $\bm{\theta}^*=\underset{\bm{\theta}}{\textrm{argmin}}\, L(\bm{\theta})$
maximizes the likelihood of the measurements.

For validation of the results, $10$ extra trajectories were generated without noise.
We now study how well it is possible to estimate the inertia and drag given this
setup using three different uncertainty models with decreasing amount of a priori knowledge.
%


\subsubsection{Uncertainty Model 1: Parameter estimation.}\label{sec:drag1}
We assume that the complete model is known, and only the drag parameters $\bm{k}$ are unknown, as well as the inertia $\bm{I}$ and the wind $\bm{w}$. We denote the estimates with hats, i.e. the full parameter vector is $\bm{\theta}=(\bm{\hat{k}},\bm{\hat{w}},\bm{\hat{I}})$.

\subsubsection{Uncertainty Model 2: Known drag structure.}\label{sec:drag2}

We now assume that we know that there is a constant wind in the global frame but that the drag model is
an unknown function of the relative speed in the air. We therefore estimate the drag with
\[
\hat{d}(\bm{v}) = \hat{n}_{\bm{\theta}_1}(\bm{v}-R^T\bm{\hat{w}})
\]
where $\hat{n}_{\bm{\theta}_1}:\mathbb{R}^3\rightarrow\mathbb{R}^3$,
to approximate the real drag
\[
d(\bm{v}) = n(\bm{v}-R^T\bm{w}),
\]
where $\bm{v}$ are the velocities in the quad-rotors frame of reference, and $\hat{n}_{\bm{\theta}_1}$ is a neural network defined by the parameters $\bm{\theta}_1$, i.e. the full set of parameters is $\bm{\theta}=(\bm{\hat{w}},\bm{\hat{I}},\bm{\theta}_1)$.

\subsubsection{Uncertainty Model 3: Black box drag.}\label{sec:drag3}

For the last case, we assume that we only know that the drag is dependent on the velocities $\bm{v}$ and angles $\bm{\eta}$, and estimate a neural network $\hat{n}_{\bm{\theta}_2}(\bm{v},\bm{\eta}):\mathbb{R}^6\rightarrow\mathbb{R}^3$. The full set of parameters is therefore $\bm{\theta}=(\bm{\hat{I}}, \bm{\theta}_2)$.
For full details of the structures of the neural networks, see Appendix \ref{app:network}.

\subsection{Solver}
To solve the differential equation we use the Tsit5 4/5 Runge-Kutta method \cite{TSITOURAS201}.
For training we use the stochastic-gradient based ADAM algorithm,
with the unusually small parameters $\beta_1 = 0.80, \beta_2=0.92$
to reflect the small size of the training set.
The learning rate is set to deceasing values in the range $(0.1, 0.0001)$.

The total number of iterations (simulations) needed for convergence is approximately 1000 for Model 3.

\subsection{Numerical Results}

We trained 10 different models for each of the uncertainty models to get an estimate of the variance of the estimation error.
Each of the 10 models were trained on the same data, but with different initial guesses for the neural networks.
We compare the root mean square (RMS) of the prediction errors,
i.e the square root of the cost function in equation \eqref{eq:cost-prediction},
on both the noisy training data $\bm{y}[k]$,
the training data without noise $\bm{p}(t_k)$,
as well as on the test data.

To estimate the approximation accuracy for the drag, the following relative error over all trajectories is used:
\begin{align}\label{eq:cost-drag}
    \sqrt{\frac{\sum_{i=1}^M\sum_{k=1}^N \|\hat{d}(\bm{v}_i(t_k))-d(\bm{v}_i(t_k))\|^2_2}{\sum_{i=1}^M\sum_{k=1}^N \|d(\bm{v}_i(t_k))\|^2_2}}
\end{align}
where $\hat{d}$ is the estimated drag model, and $\bm{v}_i$ are the true velocities from each trajectory $i$.
We evaluate this error on both the training and test data.

These values, together with the relative errors of the inertia estimates $\bm{\hat{I}}$, are presented in Table \ref{table:results}.

\subsubsection{Uncertainty Model 1:}

In all cases, the difference between the worst and best approximation is small, indicating that the method is fairly robust.
We see that the trained models are able to fit the training data
close to the expected RMSE of $0.03$ on the noisy training data.
The prediction error on both the test and training data is less than
$1\si{\milli\meter}$ on average.
The drag term is also approximated well with approximately $1\%$ relative error.
The inertia parameters $I_1,I_2$, corresponding to pitch and roll, was estimated to within approximately $2\%$,
but the inertia $I_3$ for the yaw was approximately $13\%$ from the true value.
This indicates that the simulation could be lacking sufficient excitation in these directions.
It was noted that the derivative of the cost with respect to the inertia parameter $I_3$ was an order or magnitude smaller for $I_1$ and $I_2$ when close to the true values.
It is possible that the controller is too fast around this axis to allow for proper identification.

The wind $\bm{w}$ was correctly estimated in each direction to between $0.03\%$ and $0.8\%$ with a mean of $0.4\%$ relative error.
The drag parameters $k_1,k_2$ were estimated to approximately $0.2\%$ relative error, and $k_3$ to $1.4\%$.

We now use the prediction and drag errors for this model as a base-line for the more complicated models below.

\begin{figure}
\includegraphics[width=1.0\linewidth]{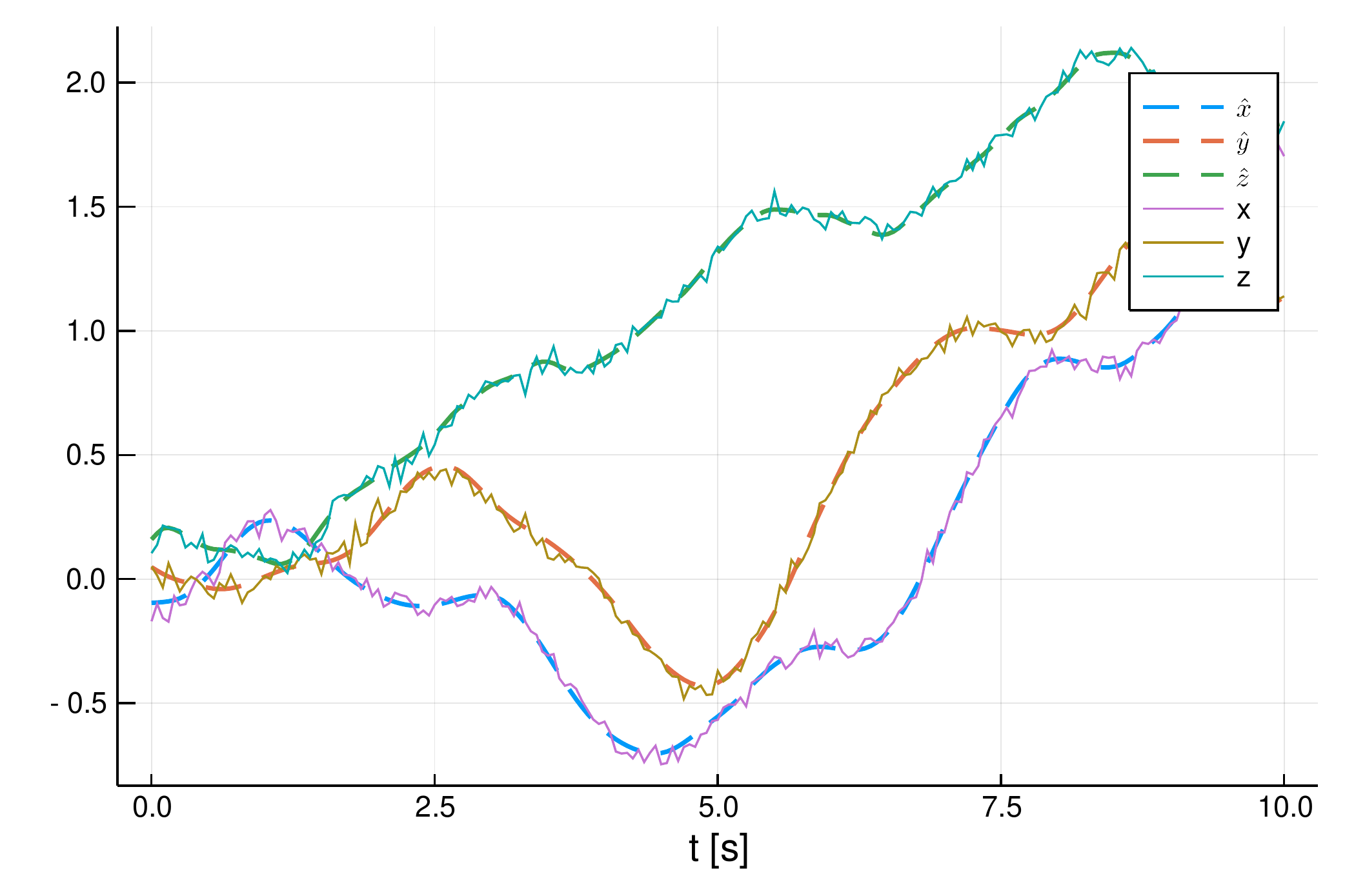}
\caption{\label{fig:est-1} Example trajectory from uncertainty model 1. Measured positions in dashed lines and predicted in solid.}
\end{figure}

\begin{figure}
\includegraphics[width=1.0\linewidth]{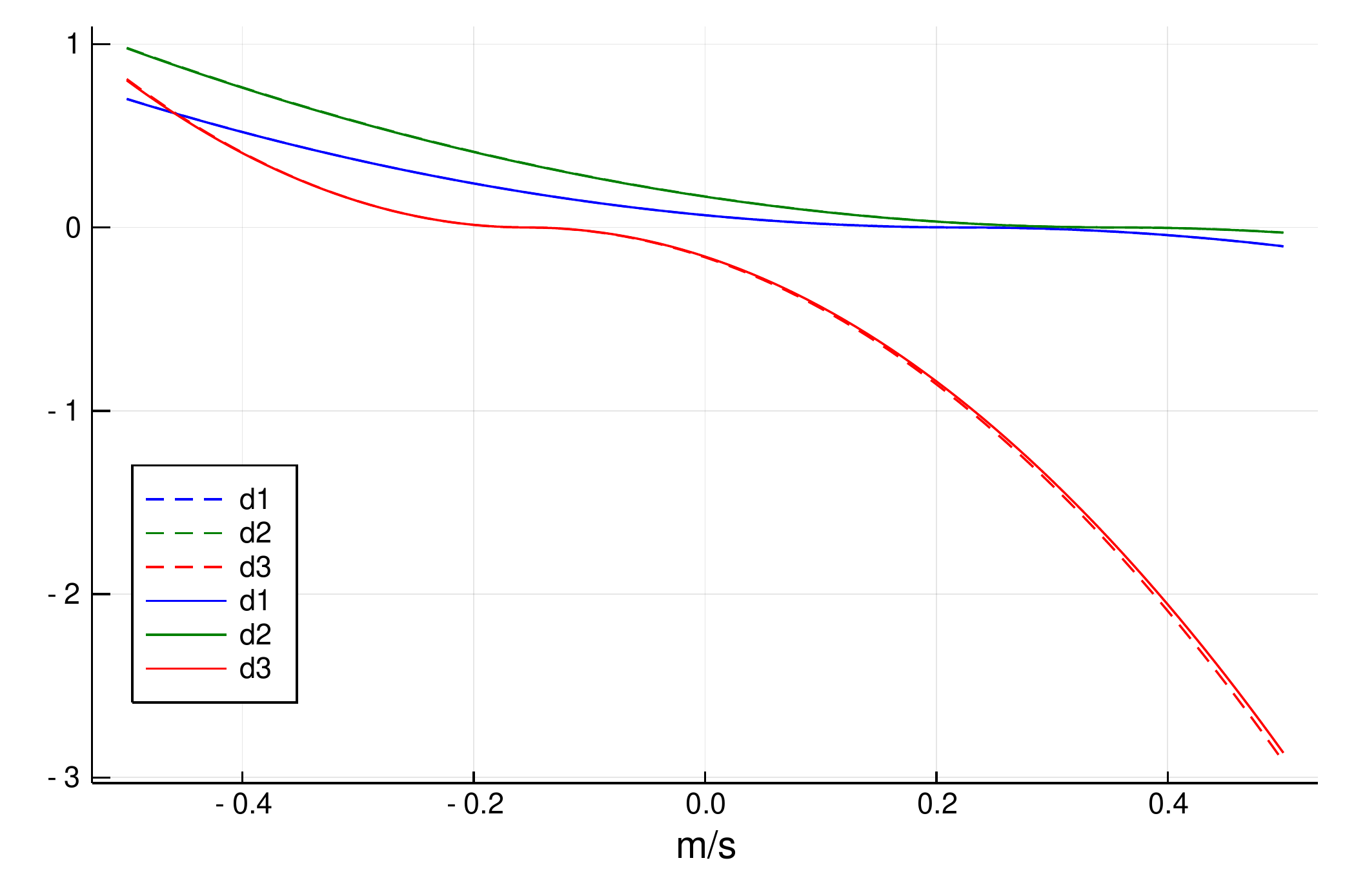}
\caption{\label{fig:drag-0} Diagonal drag terms from model in Section \ref{sec:drag1}. Estimated drag in dashed lines and real drag in solid lines}
\end{figure}

The true and predicted drag, $d(\bm{v})$ and $\hat{d}(\bm{v})$ respectively, are shown in Figure \ref{fig:drag-0} for $\bm{w}=\bm{0}$.
An example trajectory of the noisy measurements and the predicted trajectory is shown in Figure \ref{fig:est-1}.
\subsubsection{Uncertainty Model 2:}

As seen in Table \ref{table:results}, the prediction errors for this model are roughly 6 times worse than for the baseline model.
However, the errors are still considerably lower than the noise level, at approximately $0.5\si{\centi\metre}$.
The big difference compared to the baseline is the estimation of the drag term.
The error is now approximately $11\%$,
reflecting the increased difficulty of training a neural network.
The estimation error of the inertia is slightly worse in most cases with maximum errors of $2.7\%, 5.7\%$ and $10\%$ for $I_1,I_2$ and $I_3$ respectively.
The variance of the estimation errors are slightly larger,
this is probably a consequence of the initial values in the neural network.

\begin{figure}
\includegraphics[width=1.0\linewidth]{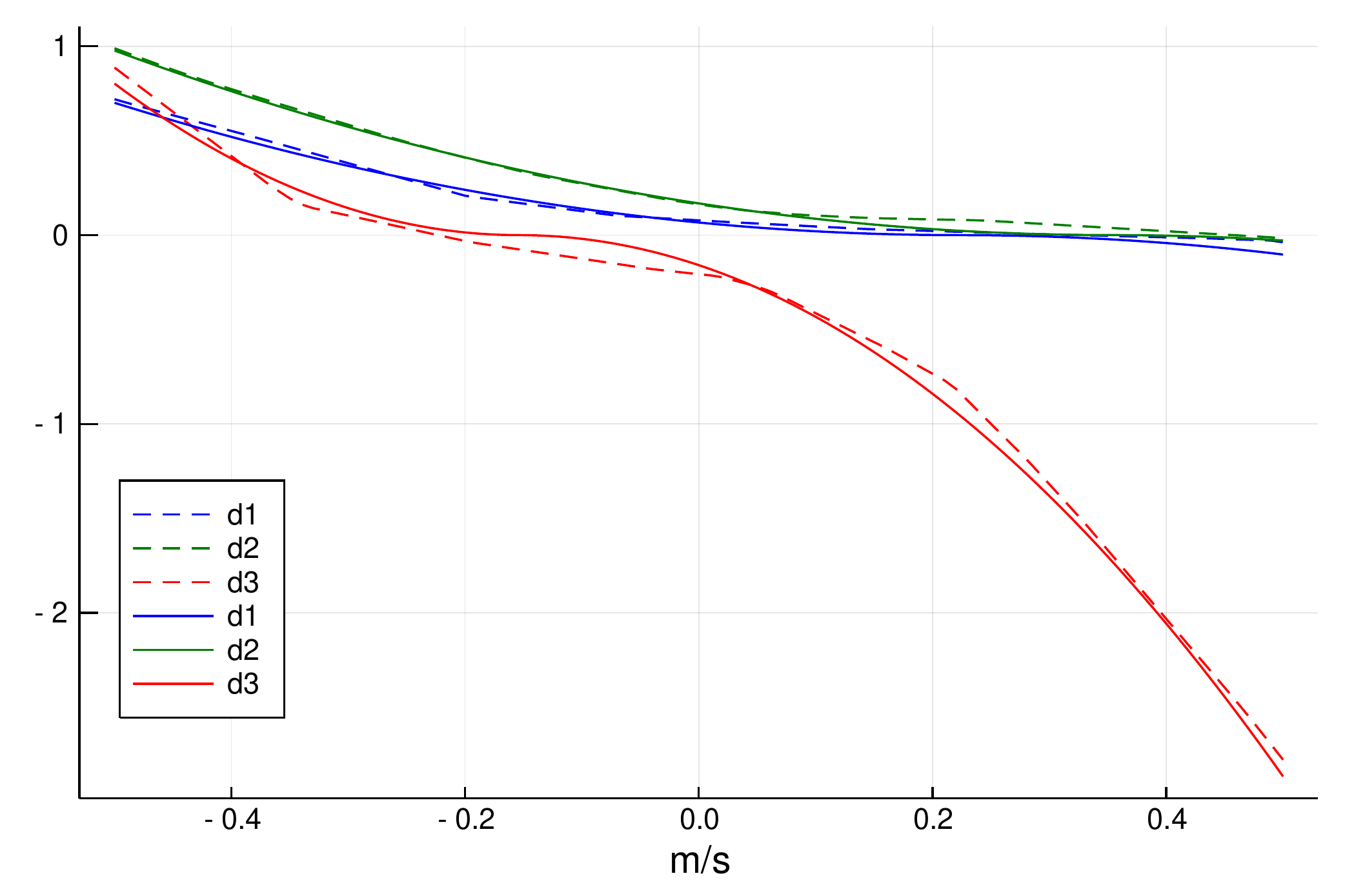}
\caption{\label{fig:drag-1} Diagonal drag terms from model in Section \ref{sec:drag2} with $R(\bm{0})=I$. Estimated drag in dashed lines and real drag in solid lines}
\end{figure}

No estimation of the wind $\bm{w}$ can be properly extracted
independently of the neural network since the network was not forced to have the property $\hat{n}_{\bm{\theta}_1}(\bm{0})=\bm{0}$.
The estimated drag is shown in Figure \ref{fig:drag-1} for $R^T(\bm{0})=I$, where $I$ is the identity matrix.

\subsubsection{Uncertainty Model 3:}

The prediction errors for this model are almost as good as for Model 2 on the training data.
However, it should be noted that the error on the training data is now lower than the noise floor ($2.95$ vs $3.0$).
This shows that we are now over-fitting to the training data,
a consequence of the more expressive model.
As a result, the prediction error on the test data, is on average a factor 2 worse than model 2.
We also see that the estimated drag function does not always
generalize well to the test data, with a maximum relative error of $45\%$.
The variance of the inertia estimates is now much larger,
and in some cases even match the true values better than model 1.

\begin{figure}
\includegraphics[width=1.0\linewidth]{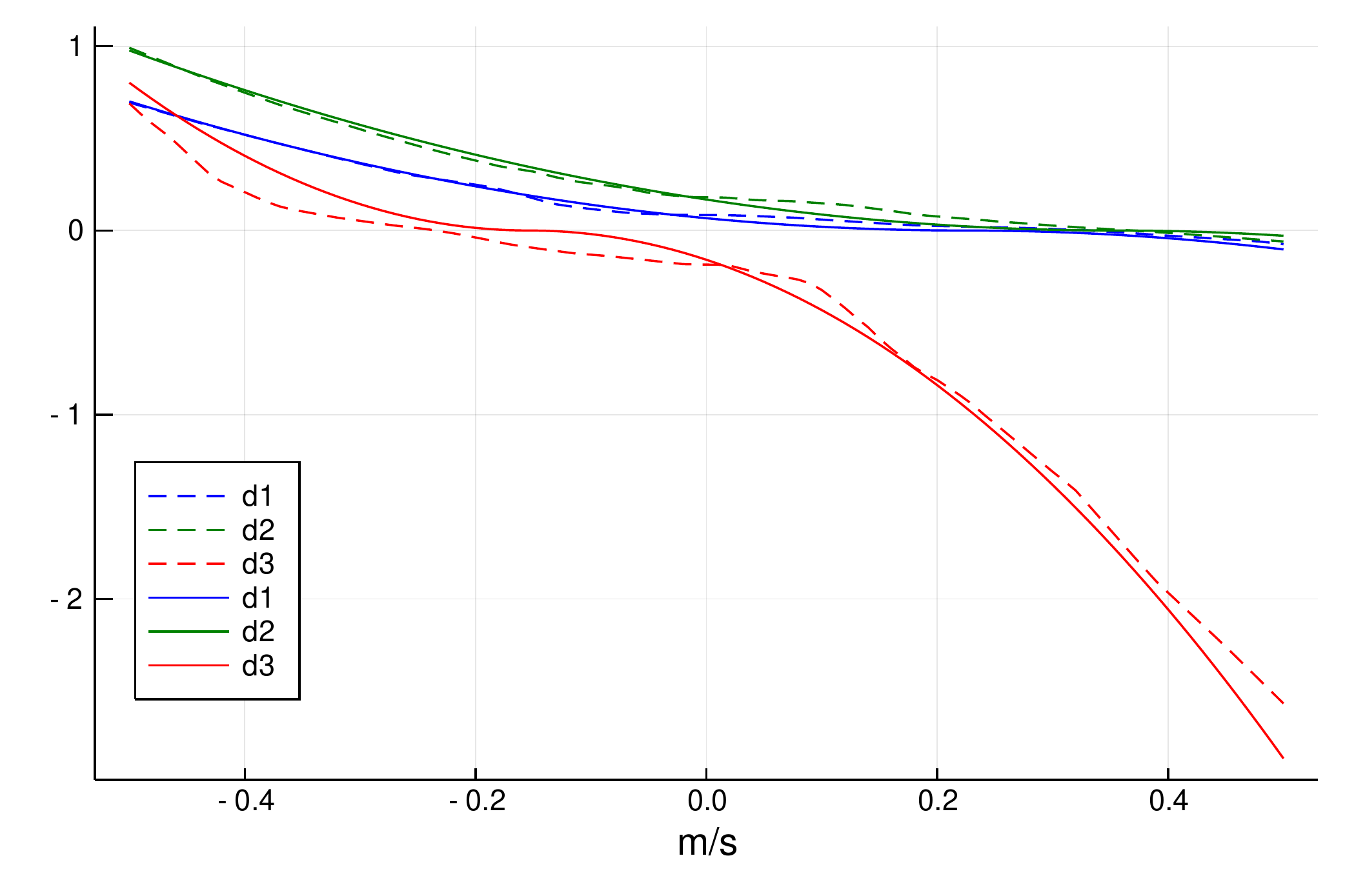}
\caption{\label{fig:drag-2} Diagonal drag terms from model in Section \ref{sec:drag3} with $\bm{\eta}=\bm{0}$. Estimated drag in dashed lines and real drag in solid lines}
\end{figure}

The estimated drag function is shown in Figure \ref{fig:drag-2}.

%

%
%

\begin{table*}[!h]
\centering
\hspace*{-0.2\textwidth}
\begin{tabular}{|c||c|c|c|}
\hline
\diagbox{Variable}{Model} & $- k_i\cdot v_i|v_i|$ &   $\hat{n}(\bm{v}-R^T\hat{\bm{w}})$ & $\hat{n}(\bm{m},\bm{\eta})$ \tabularnewline
\hline
\hline
Train Error RMSE (noisy)& $(3.01,3.01,3.01)\cdot10^{-2}$ & $(3.00,3.00,3.00)\cdot10^{-2}$  & $(2.95,2.96,2.97)\cdot10^{-2}$  \tabularnewline
\hline
Train Error RMSE (no noise) & $(6.71,7.21,7.53)\cdot10^{-4}$ & $(3.74,3.92,4.04)\cdot10^{-3}$ & $(5.19,5.41,5.80)\cdot10^{-3}$  \tabularnewline
\hline
Test Error RMSE & $(8.02,8.64,9.10)\cdot10^{-4}$ & $(5.19,5.71,6.18)\cdot10^{-3}$ & $(5.80,11.3,19.1)\cdot10^{-3}$  \tabularnewline
\hline
Drag train relative error & $(1.26,1.32,1.37)\cdot10^{-2}$ & $(0.106,0.112,0.118)$ & $(0.15,0.19,0.27)$ \tabularnewline
\hline
Drag test relative error & $(1.27,1.33,1.38)\cdot10^{-2}$ & $(0.105,0.111,0.116)$ & $(0.16,0.23,0.45)$ \tabularnewline
\hline
$I_{1}$ & $(0.0192,0.0196,0.0207)$ & $(0.00198,0.0109,0.0274)$ & $(0.001,0.013,0.030)$ \tabularnewline
\hline
$I_{2}$ & $(0.0102,0.0149,0.0178)$ & $(0.0214,0.0414,0.0571)$ & $(0.012,0.046,0.081)$ \tabularnewline
\hline
$I_{3}$ & $(0.129,0.136,0.143)$ & $(0.120,0.153,0.194)$ & $(0.013,0.094,0.222)$ \tabularnewline
\hline
\end{tabular}
\caption{\label{table:results}Errors for the three different uncertainty models are presented above.
10 different models were identified for each case and the $(\min, \textrm{mean}, \max)$ of the errors are shown.}
\end{table*}
\section{Conclusions}
We have shown how to perform continuous-time system
identification using neural networks for hybrid systems.
We studied a quad-rotor example with the goal of estimating  the drag.
We show that it is possible to get a good model for the drag and relatively good estimates of the inertia parameters with the right prior knowledge.
As expected, the prediction error increases slightly with decreasing model knowledge,
but even in the case of very limited knowledge, it is still possible to get a decent drag model for prediction.
For future research it would be interesting to see how these methods perform on real-world data, where noise has to be considered in the state-feedback.
\section*{Acknowledgement}
The authors thanks Sebastian Banert and Marcus Greiff for valuable comments and insights.

\appendix
\section{Data Generation}\label{app:data}
The initial state $\bm{x}(0)$, reference signals $\bm{r}[k]$, input noise $\bm{r_n}[k]$, and measurement noise $\bm{n}[k]$ was generated according to the following distributions
\begin{align*}
(\bm{x}(0))_i &\sim \mathcal{N}(0,\,\sigma_i^{2}), \textrm{ where } \sigma_i = \begin{cases}
 0.02 &\textrm{ for } i\in\{1,2,3\}\\
 0.2  &\textrm{ for } i\in\{4,\dots,12\}
\end{cases}\\
(\bm{r}[k])_i &= \sum_{j=1,\dots,k} r_{i,j}, \textrm{ where } r_{i,j} \sim \mathcal{N}(0,\,0.1^{2})\\
(R_d[k])_i &\sim \mathcal{N}(0,\,10^{2})\\
(\bm{n}[k])_i &\sim \mathcal{N}(0,\,0.03^{2}).
\end{align*}
The initial guesses, where applicable, were
$\bm{\hat{w}}=(0,0,0)$,
$\bm{\hat{k}} = (0.027,0.027,0.162)$ and $\bm{\hat{I}}=(6.0, 7.0, 11.0)$.
\section{Network Structure}\label{app:network}
For both uncertainty model 2 and 3, networks with fully connected dense layers
with LeakyReLu activation functions are used, without any activation function on the output layer.
For model 2, a wide network structure with a single hidden layer with 120 nodes is used to simplify the training.
Model 3 has 4 hidden layers with widths (20, 120, 30) in the order from input to output layer.

\bibliographystyle{agsm}
\bibliography{root}

\end{document}